\input amstex 
 \documentstyle{amsppt}
\NoRunningHeads 
 
\magnification 1200
%
%
\catcode`@=11
\def\binrel@#1{\setbox\z@\hbox{\thinmuskip0mu
\medmuskip\m@ne mu\thickmuskip\@ne mu$#1\m@th$}%
\setbox\@ne\hbox{\thinmuskip0mu\medmuskip\m@ne mu\thickmuskip
\@ne mu${}#1{}\m@th$}%
\setbox\tw@\hbox{\hskip\wd\@ne\hskip-\wd\z@}}
\def\overset#1\to#2{\binrel@{#2}\ifdim\wd\tw@<\z@
\mathbin{\mathop{\kern\z@#2}\limits^{#1}}\else\ifdim\wd\tw@>\z@
\mathrel{\mathop{\kern\z@#2}\limits^{#1}}\else
{\mathop{\kern\z@#2}\limits^{#1}}{}\fi\fi}
\def\underset#1\to#2{\binrel@{#2}\ifdim\wd\tw@<\z@
\mathbin{\mathop{\kern\z@#2}\limits_{#1}}\else\ifdim\wd\tw@>\z@
\mathrel{\mathop{\kern\z@#2}\limits_{#1}}\else
{\mathop{\kern\z@#2}\limits_{#1}}{}\fi\fi}
\def\circle#1{\leavevmode\setbox0=\hbox{h}\dimen@=\ht0
\advance\dimen@ by-1ex\rlap{\raise1.5\dimen@\hbox{\char'27}}#1}
\def\sqr#1#2{{\vcenter{\hrule height.#2pt
     \hbox{\vrule width.#2pt height#1pt \kern#1pt
       \vrule width.#2pt}
     \hrule height.#2pt}}}
\def\square{\mathchoice\sqr34\sqr34\sqr{2.1}3\sqr{1.5}3}
\def\force{\hbox{$\|\hskip-2pt\hbox{--}$\hskip2pt}}  
\catcode`@=\active

%
%
\baselineskip 18pt
\define\un{\underbar}
\define\var{\vartheta}
\define\pmf{\par\medpagebreak\flushpar}
\topmatter
\title
A note on canonical functions
\endtitle
\author
Thomas Jech and Saharon Shelah \endauthor
\thanks Supported by NSF and by a Fulbright grant; and
Publ. 378, partially supported by the B.S.F.
\endthanks
\endtopmatter

\centerline{Abstract}
\par
We construct a generic extension in which the $ \aleph_2$ nd canonical function on
$ \aleph_1$ exists.
\par
\newpage

\par \flushpar
\un{Introduction}  For ordinal functions on $ \omega_1$, let $f < g$ if 
$ \{ \xi < \omega_1: f ( \xi ) < g ( \xi ) \}$  contains a closed unbounded
set.  By induction on  $ \alpha,$ the $ \underline{\alpha}$\un{th canonical function}
$f_{ \alpha}$ is defined (if it exists) as the least ordinal function 
greater than each $f_{ \beta}, \beta < \alpha$ (i.e. if $h$ is any other function greater than all $f_\beta, \beta < \alpha$, then $f_\alpha \le h$).
If $f_{ \alpha} $
exists then it is unique up to the equivalence  
$$ \{ \xi < \omega_1 : f ( \xi ) = g ( \xi ) \}  \; \text{contains a closed unbounded
set.} $$  
It is well known [2] that for each  $ \alpha < \omega_2,$ the 
$ \alpha$th canonical function exists. A. Hajnal has shown (private communication)
that if $V = L$ then the $ \aleph_2$nd canonical function does not exist.
In this note we show that it is consistent that the $ \aleph_2$nd canonical 
function exists. We prove a somewhat more general result:
\pmf
\un{THEOREM}  Assume that $2^{ \aleph_0} = \aleph_1,$ and let $ \var$ be an
ordinal.  There is a cardinal preserving generic extension in which for
every $ \alpha < \var$ the $ \alpha$th canonical function exists.
\pmf
\un{Remarks}
\par
\flushpar
1.  In the model of the theorem, each $f_{\alpha} , \alpha < \var,$ is a function 
into $ \omega_1;$ thus $2^{ \aleph_1} \ge | \var |$
\par \flushpar
2.  In the model of [3], canonical functions exist for all ordinals $ \alpha.$
The model is constructed under the assumption of a measurable cardinal;
that assumption is necessary since if all canonical functions exist then the 
closed unbounded filter is precipitous.
\par \flushpar
3. Consider the statement
\flushpar
``the constant function $ \omega_1$ is a canonical function"
\flushpar
Its consistency implies the consistency of set theory with predication [4], 
and hence of various mildly large cardinals.
\pmf
4. The theorem generalizes, in the obvious way, to ordinal functions on any 
regular uncountable cardinal.  
\par \medpagebreak \par 
The  proof of the theorem uses iterated forcing.  We use the standard 
terminology of forcing; see e.g. [1] for iterated forcing.  A notion of
forcing is $ \omega $\un{-distributive} if it adds no new countable sequences of
ordinals; it is $ \aleph_2$-c.c. if it has no antichain of size
$ \aleph_2.$    A set $S \subseteq \omega,$ is \un{costationary} if 
$ \omega_1 - S$ is stationary. 
 
By a \un{countable model} $N$ we mean a countable  elementary submodel of
$(V_{ \kappa}, \in )$  where $ \kappa$ is a sufficiently large cardinal.  A sequence $ \{
p_n \}_{n  \in \omega}$  of conditions in $P$ is \un{generic} for  a countable model
$N$ if $ P \in N, \; \; \{ p_n : n \in \omega \} \subset N,$ and if
$ \{ p_n \}_n $ meets every dense set $D \subseteq P$ such that $D \in N.$
\pmf
\un{Construction of the model}
 
We construct the forcing $P$ in two stages:  first we adjoin generically
a $ \var$-sequence of functions $f_i :  \omega_1 \to \omega_1, \; \quad i < \var,$
such that $f_i < f_j $ whenever $i < j.$  The forcing $P_0 $
that does it is $ \omega$-closed and satisfies the $ \aleph_2$-chain condition.
The second stage is an iteration, with countable support, of length $ \lambda =
(2^{ \aleph_1} \cdot | \var |)^+$ that successively destroys all stationary sets which
witness that the functions $f_i$ are not canonical.  We will prove that the
iteration forcing is $ \omega$-distributive and $ \aleph_2$-c.c.  Hence $P$
preserves cardinals, and one can arrange all the names for subsets of $ \omega_1$
in a sequence $ \{ \dot S_{ \alpha} : 1 \le \alpha < \lambda \} $ 
such that for each $ \alpha, \; \dot S_{ \alpha}$ is in $ M_{ \alpha} = V^{P 
|  \alpha}.$  Moreover, this can be done in such a way that each
$ \dot S$ appears in the sequence cofinally often.  We remark that if 
$ S \subseteq \omega_1$ is in $M_{\lambda}$ then $S \in M_{ \alpha}$ for some $ \alpha < \lambda;$
if $M_{ \lambda} \models S$ is stationary then $M_{ \alpha} \models S$
is stationary; if $M_{\alpha} \models f < g$ then $M_{\lambda} \models f < g,$ 
and if $M_{ \lambda} \models f < g$ then for all sufficiently large
$ \alpha < \lambda, \; M_{ \alpha} \models f < g.$
 
\par \flushpar
\un{Definition of} $P_0$
 
A condition consists of
\item{(a)} a countable ordinal $ \gamma$
\item{(b)} a countable set $A \subset \var$
\item{(c)} closed subsets $c_{ij} \; \text{of} \; \gamma \; \; (i, j \in A, \; i < j)$
\item{(d)} functions $f_i: \; \gamma \rightarrow \omega_1,$
\par \flushpar
such that for all $ i, j  \in A, \; i < j, \; f_i ( \xi) < f_j
( \xi)$ for all $ \xi \in c_{i j}$
\pmf
A stronger condition increases $ \gamma$ and $A,$  extends the $f_i,$
and end-extends the $c_{ij}.$  
 
The forcing $P_0$ is $ \omega$-closed,
and is $ \aleph_2-c.c.$ because $2^{ \aleph_0}  = \aleph_1.$
Let $ \dot f_i,\; i < \var,$ denote the names for the generic functions forced
by $P_0.$  Clearly, $M_0 \models \dot f_i < \dot f_j$
whenever $ i < j.$
\pmf
\un{Definition of P}
 
\par \flushpar
$ P = P_{ \lambda} $ is an iteration with countable support.  For
$ 1 \le \alpha \le \lambda, \;  P_{ \alpha}$ is the set of all $ \alpha-$sequences
 $ \{ p( \beta):  \beta < \alpha \}$ with countable support such that $p (0) \in
 P_0$ and such that $p( \beta) = \emptyset$
( trivial condition) unless the following is forced by $p | \beta:$
\item{(1)} for some $ i < \var$ there exists a function $g$ such that 
$g > \dot f_j$ for all $j < i,$ and $g ( \xi) < \dot f_i ( \xi)
$ everywhere on $ \dot S_{ \beta}.$  
\par \flushpar
In that case $p( \beta) $ is a countable
closed set of countable ordinals that is forced by $p | \beta $ to be disjoint
from $ \dot S_{ \beta}.$
 
A condition $q$ is stronger than $p$ if $q (0) \le p (0)$ and for all $ \beta, \; 1 \le \beta < \alpha, \;
q( \beta)$ end-extends $p ( \beta).$  
\par \flushpar
For every $ \beta $ that satisfies (1),
the forcing produces a closed subset $C$ of $ \omega_1 $ disjoint from $ \dot S_\beta,$ 
and $C$ is unbounded as long as $ ( \omega_1 - \dot S_{ \beta}) $ is unbounded.
We shall prove that $P_{ \lambda}$ is $ \omega-$distributive and $ \aleph_2-c.c.,$
and that in $M_{ \lambda}$ the functions $ \dot f_i$ are canonical.
\pmf
\un{Lemma}  Let $N$ be a countable model such that $ P \in N,$ and let
$ \delta = \omega_1 \cap N.$  If $ \{ p_n \}_{n \in \omega}$ is a generic sequence 
 for $N,$ then there exists a $q$ stronger than all the $p_n,$ and such that for all 
$i \in N, \; q$ forces
\item{(2)}    $ \qquad \qquad  \dot f_i ( \delta) = sup \{ \dot f_j
( \delta) + 1: \; \; j < i \; \text{and} \; j \in N \}$
\pmf
Proof.  Let $X$ be the union of the supports of $p_n, \; n \in \omega;$ 
note that $X \subset N.$  We construct $q ( \beta) $ by induction on $  \beta.$
If $ \beta \notin X$ we let $q( \beta) = \emptyset.$
 
First let $ \beta = 0.$  Look at $ \{ p_n (0) \}_{  n \in \omega}.$  By the
genericity  of the sequence, the ordinals $ \gamma_n$ converge to $ \delta,$
 the countable sets $A_n$ converge to $A = \var  \cap N,$ the closed sets 
$(c_{ij})_n $ converge to $c_{ij} \subseteq \delta$ and the functions $(f_i)_n$
converge to functions $ f_i: \; \delta \to \delta$ such that $ f_i < f_j$ on
$c_{ij}.$  
\par \flushpar
Let $ \gamma = \delta + 1,$ let $ \bar c_{ij} = c_{ij} \cup \{ \delta \}$
and let $ \bar f_i$ be the extensions of the $ f_i$s that satisfy (2).  Let
$q(0)$ be the condition $ ( \delta, A, \bar c_{ij},  \bar f_i);$ clearly, $q(0)$
forces (2).  
\par \flushpar
Now let $1 \le \beta < \lambda, \; \beta \in X,$ and assume that we 
have already constructed $ q | \beta$ stronger than all the $p_n | \beta.$  As
eventually all $ p_n | \beta$ force (1), it follows by their genericity that
the countable sets $ p_n ( \beta)$ converge to a closed subset of $ \delta.$
So we let $q (\beta) = \underset n \to \bigcup p_n ( \beta) \cup \{ \delta \},$
and in order that $q$ be a condition, we have to verify that 
$q | \beta \force \delta \notin \dot S_{ \beta}.$  
\par \flushpar
Let $q'$ be any
condition in $P_{ \beta}$  stronger than $q |  \beta.$  Since
$ \beta \in N,$ we may assume that $ \dot S_{ \beta} \in N,$ and $N$ satisfies 
that for eventually all $n, p_n | \beta$ forces (1).  It follows that there exists  a
condition $ r \le q',$ some $ i \in N$ and some $ \dot g \in N \cap M_{ \beta}$
such that for all $ j < i$ in $N,\; \; r \force \dot f_j < 
\dot g,$ and $ r \force ( \forall \xi \in \dot S_{ \beta} ) \dot g ( \xi) 
< \dot f_i ( \xi).$
\par \flushpar
For each $j < i $ in $N,$ there exists an $M_{ \beta}$-name
$ \dot C_j \in N$ such that every $p \in P_{ \beta} $ forces that
$ \dot C_j$ is closed unbounded, and 
$ r \force ( \forall \xi \in \dot C_j) \dot f_j ( \xi) < \dot g ( \xi).$ 
It follows, by the genericity of  
$ \{ p_n \}_n,$ that $q \force \dot C_j \cap \delta $ is cofinal in $ \delta,$ and
 so $ q \force \delta \in \dot C_j $  Hence $r$ forces that for all
$j < i$ in $N, \; \; \dot f_j ( \delta) < \dot g ( \delta).$  
But since $r$ also forces (2), it forces $ \dot f_i ( \delta) \le \dot g ( \delta),
$ and therefore $ r \force \delta \notin \dot S_{ \beta.}$  [Note that the 
proof also yields that $q | \beta$ forces that $ \dot S_{ \beta} $
is costationary, as the argument above proves that $ q | \beta \force \delta \in \dot C$
for every club name in $N$.] $ \square$ 
\pmf
\un{Corollary} $P$ is $ \omega$-distributive and $ \aleph_2-$c.c.
\par \flushpar
Proof  If $ \dot X \in M_{ \lambda}$ is a name for a countable set of ordinals
and $ p \in P_{ \lambda},$ let $N$ be a countable model such that $ \dot X
\in N, \quad P_{ \lambda} \in N$ and $ p \in N.$   Let $ \{ p_n \}_n$ be a
generic sequence for $N$ such that $ p_0 = p.$  By Lemma $ \{ p_n \}_n$
has a lower bound $q,$ and by genericity, $q$ decides each $ \dot X (n).$
Hence $P_{ \lambda}$ is $ \omega-$distributive.
 
For each $ \alpha, M_{ \alpha + 1}$ is a forcing extension of $M_{ \alpha}$
via a set of conditions of size $ \aleph_1,$ therefore $ \aleph_2-$c.c.  As
each $P_{ \alpha}$ is an iterated forcing with countable support, it satisfies
the $ \aleph_2-$c.c. as well.  \hfill $ \square$
 
We shall finish the proof of the Theorem by showing that in the generic extension
by $P,$ the functions $f_i, i < \var,$ are canonical.  We show that for each
$ i < \var, \; \; f_i $ is the least function greater than all the
$f_j, \; \; j < i.$  We already know that $f_i > f_j$ for all $j < i.$
 
Let $g \in M_{ \lambda}$ be any function such that $f_j < g$ for all $j < i,$ and
let $S = \{ \xi : \; g ( \xi ) < f_i ( \xi ) \}.$  We want to show that $S$
is nonstationary.  Let $ \beta$ be an ordinal such that $S_{ \beta} = S,$
sufficiently large so that all the clubs witnessing $f_j < g$  (all $j<i$)
belong to $M_{ \beta}.$  Hence $M_{ \beta}$ satisfies (1), and so the
forcing at stage $ \beta$ adjoins a closed unbounded  set that is disjoint from
$S$. 
\par
\midspace{1cm}
\flushpar
\un{Acknowledgment} The first author appreciates the hospitality of the Hebrew
Univeristy Mathematics Department during his sabbatical leave.
\par
\newpage
\par \flushpar
\un{References}
\pmf
\item {[1]} J. Baumgartner, Iterated forcing, in:  Surveys in set theory,
London Math. Soc. Lecture Note Ser. \un{87} (1983), p. 1 - 59.
\item {[2]} F. Galvin and A. Hajnal, Inequalities for cardinal powers,
Annals of Math. \un{101} (1975), 491-498.
\item{[3]} T. Jech, M. Magidor, W. Mitchell and K. Prikry, Precipitous Ideals,
J. Symb. Logic \un{45} (1980), 1-8
\item {[4]} T. Jech and W. Powell, Standard models of set theory with predication,
Bull. Amer. Math. Soc. \un{77} (1971), p. 808-813
\par
\midspace{1cm}
\par \flushpar
The Pennsylvania State University    \hfill  The Hebrew University in Jerusalem
 
\bye